\documentclass[10pt,letterpaper]{article}
\usepackage{blindtext,titlefoot}
\usepackage{authblk}
\usepackage{marvosym}
\usepackage[utf8]{inputenc}
\usepackage[english]{babel}
\usepackage{amsmath}
\usepackage{amsfonts}
\usepackage{amssymb}
\usepackage{makeidx}
\usepackage{graphicx}
\usepackage{lmodern}
\usepackage{kpfonts}
\usepackage{dsfont}
\usepackage{cancel}
\usepackage{amsthm} 
\usepackage[left=2cm,right=2cm,top=2cm,bottom=2cm]{geometry}
\usepackage{subcaption}
\usepackage{multirow}
\usepackage{float}
\usepackage{url}
\usepackage{breakurl}
\usepackage[breaklinks]{hyperref}
\usepackage{multicol}
\usepackage{mathtools, cuted}
\usepackage{lipsum}
\usepackage{booktabs}
\usepackage{comment}
\usepackage{cite}

\providecommand{\nset}[1]{
\mathbb{#1}
}
\providecommand{\set}[1]{
\left\{#1\right\}
}

\providecommand{\ifr}[5]{
{}^{#1}_{#2}{#3}_{#4}^{#5}
}
\providecommand{\gam}[1]{
\Gamma\left(#1 \right)
}
\providecommand{\re}[1]{
\hbox{Re}\left(#1 \right)
}
\providecommand{\im}[1]{
\hbox{Im}\left(#1 \right)
}
\providecommand{\norm}[1]{
\left\lVert #1 \right\rVert
}
\providecommand{\abs}[1]{
\left\lvert #1 \right\rvert
}
\providecommand{\ds}[1]{
\displaystyle #1
}
\providecommand{\der}[3]{
\dfrac{#1^{#3} }{ #1 #2^{#3}}
}

\providecommand{\kr}[1]{
\hbox{ker}\left(#1\right)
}

\providecommand{\rnd}[2]{
\hbox{Rnd}_#2\left(#1\right)
}

\newtheorem{theorem}{ Theorem}[section]
\newtheorem{definition}[theorem]{Definition}
\newtheorem{proposition}[theorem]{Proposition}
\newtheorem{corollary}[theorem]{Corollary}
\newtheorem{example}[theorem]{Example}

\usepackage{enumitem} 
\setlist[itemize]{noitemsep} 

\usepackage{titlesec} 
\titleformat{\section}[block]{\large\bfseries\scshape\centering}{\thesection.}{1em}{} 
\titleformat{\subsection}[block]{\large\bfseries\scshape\centering}{\thesubsection.}{1em}{}
\titleformat{\subsubsection}[block]{\large\bfseries\scshape\centering}{\thesubsubsection.}{1em}{} 

\title{\huge\bfseries An approximation to zeros of the Riemann zeta function using fractional calculus}

\author[,a]{A. Torres-Hernandez  \footnote{Email: anthony.torres@ciencias.unam.mx; ORCID: 0000-0001-6496-9505}}
\affil[a]{Department of Physics, Faculty of Science - UNAM, Mexico}

\author[,b]{F. Brambila-Paz \footnote{Email: fernandobrambila@gmail.com; ORCID: 0000-0001-7896-6460}}
\affil[b]{Department of Mathematics, Faculty of Science - UNAM, Mexico}

\date{}

\begin{document}

\maketitle


\begin{abstract}

In this document, as far as the authors know, an approximation to the zeros of the Riemann zeta function has been obtained for the first time using only derivatives of constant functions, which was possible only because a fractional iterative method was used. This iterative method, valid for one and several variables, uses the properties of fractional calculus, in particular the fact that the fractional derivatives of constants are not always zero, to find multiple zeros of a function using a single initial condition. This partly solves the intrinsic problem of iterative methods that if we want to find $N$ zeros it is necessary to give $N$ initial conditions. Consequently, the method is suitable for approximating nontrivial zeros of the Riemann zeta function when the absolute value of its imaginary part tends to infinity.  The deduction of the iterative method is presented, some examples of its implementation, and finally $53$ different values near to the zeros of the Riemann zeta function are shown.

\end{abstract}

\section{Introduction}

A classic problem in mathematics, which is of common interest in physics and engineering, is finding the set of zeros of a function $f:\Omega \subset \nset{R}^n \to \nset{R}^n$, that is,

\begin{eqnarray}\label{eq:1-001}
\kr{f}:= \set{\xi \in \Omega \ : \ \norm{f(\xi)}=0},
\end{eqnarray}

where $\norm{ \ \cdot \ }: \nset{R}^n \to \nset{R}$ denotes any vector norm. Although finding the zeros of a function may seem like a simple problem, many times it involves solving an \textbf{algebraic equation system} given by the following expression

\begin{eqnarray}\label{eq:1-002}
\left\{
\begin{array}{c}
\left[f\right]_1(x)=0\\
\left[f\right]_2(x)=0\\
\vdots \\
\left[f\right]_n(x)=0
\end{array}\right. ,
\end{eqnarray}

where $[f]_k: \nset{R}^n \to \nset{R}$ denotes the $k$-th component of the function $f$. It should be mentioned that system \eqref{eq:1-002} may represent a \textbf{ linear system} or a \textbf{nonlinear system}. The function that has the set of zeros of most interest in mathematics,  being possibly the most famous set of zeros,  corresponds to the  \textbf{Riemann zeta function}

\begin{eqnarray*}
\zeta(x)=\sum_{k=1}^\infty \dfrac{1}{k^x}. 
\end{eqnarray*}

The Riemann hypothesis, under the assumption that it is true, is responsible for giving a general expression to the zeros of function $\zeta$, and may be written compactly as

\begin{eqnarray}\label{eq:1-004}
\begin{array}{cccccc}
\mbox{ If \hspace{0.1cm} }\xi \in \nset{C}& \Rightarrow &
\xi \in \kr{\zeta} & \Leftrightarrow & \xi=\left\{
\begin{array}{cc}
-2x, &\mbox{ for all } x\in \nset{N}\\
0.5+x i,& \mbox{ for some } x\in \nset{R}
\end{array}\right.,
\end{array}
\end{eqnarray}

when $\xi\neq -2x$ is known as a \textbf{nontrivial zero} of the function $\zeta$. In general, it is necessary to use numerical methods of the iterative type to construct a sequence  $\set{x_i}_{i=0}^\infty$, that under certain conditions, allows to approximate 
to values  $\xi\in \kr{\zeta}$ with $\xi\neq-2x$, that is,

\begin{eqnarray*}
x_i\underset{i\to \infty}{\longrightarrow} \xi \in \kr{\zeta}.
\end{eqnarray*}

If we want to find $N$ values $\xi\in \kr{\zeta}$, it is necessary to give $N$ initial conditions $x_0$. The previously described is an intrinsic problem of iterative methods, because time must first be spent determining to initial conditions $x_0$ before beginning to search the values $\xi$. Furthermore, it is sometimes necessary that the initial conditions are near to the searched values to guarantee convergence, that is,

\begin{eqnarray}
\abs{x_0-\xi}<\epsilon &\Leftrightarrow & x_i\underset{i\to \infty}{\longrightarrow} \xi \in \kr{\zeta}.
\end{eqnarray}

These problems are partially solved using \textbf{fractional iterative methods}, because they can determine $N$ values $\xi$ using a single initial condition $x_0$, and the initial condition does not need to be near to the searched values. In this document is used a fractional iterative method does not explicitly depend on the fractional derivative of the function for which zeros are sought, then it is an ideal iterative method for working with functions that can be expressed in terms of a series or for solving nonlinear systems in several variables.

\section{Fixed Point Method}

Let  $\Phi:\nset{R}^n \to \nset{R}^n$ be a function. It is possible to build a sequence $\set{x_i}_{i=0}^\infty$ defining the following iterative method

\begin{eqnarray}\label{eq:2-001}
x_{i+1}:=\Phi(x_i),
\end{eqnarray}

if it is true that the sequence $\set{x_i}_{i=0}^\infty$ converges to a value $\xi\in \nset{R}^n$, and if the functions  $\Phi$ is continuous around $\xi$, it holds that

\begin{eqnarray}\label{eq:2-002}
\xi=\lim_{i\to \infty}x_{i+1}=\lim_{i\to \infty}\Phi(x_i)=\Phi\left(\lim_{i\to \infty}x_i \right)=\Phi(\xi),
\end{eqnarray}

the above result is the reason by which the method \eqref{eq:2-001} is known as \textbf{fixed point method}. Furthermore, the function $\Phi$ is called an \textbf{iteration function}. To understand the nature of the convergence of the iteration function $\Phi$, the following definition is necessary \cite{plato2003concise}:

\begin{definition}
Let $\Phi:\nset{R}^n \to \nset{R}^n$ be an iteration function. The method \eqref{eq:2-001} to determine  $\xi\in \nset{R}^n$, is called  \textbf{(locally) convergent}, if it there exists $\delta>0$ such that for all initial condition

\begin{eqnarray*}
x_0\in B(\xi;\delta):=\set{y\in \nset{R}^n \ : \  \norm{y-\xi}<\delta},
\end{eqnarray*}

it holds that

\begin{eqnarray}\label{eq:2-003}
\lim_{i \to \infty}\norm{x_i-\xi}\to 0 & \Rightarrow & \lim_{i\to \infty}x_i=\xi.
\end{eqnarray}

\end{definition}

Before continuing, it is necessary to define the order of convergence of an iteration function $\Phi$ \cite{plato2003concise}.

\begin{definition}
Let $ \Phi: \Omega \subset \nset{R}^ n \to \nset{R}^ n $ be an iteration function with a fixed point $ \xi \in \Omega $. Then the method \eqref{eq:2-001} is called  \textbf{(locally) convergent of (at least) order $ \boldsymbol{p} $} ($ p \geq 1 $), if there are exists $ \delta> 0 $  and $ C $ a non-negative constant,  with $ C <1 $ if $ p = 1 $, such that for any initial value $ x_0 \in B (\xi; \delta) $ it holds that

\begin{eqnarray}\label{eq:c2.08}
\norm{x_{k+1}-\xi}\leq C \norm{x_k-\xi}^p, & k=0,1,2,\cdots,
\end{eqnarray}

where $ C $ is called convergence factor.

\end{definition}

The order of convergence is usually related to the speed at which the sequence generated by \eqref{eq:2-001} converges. The following theorem, allows characterizing the order of convergence of an iteration function $ \Phi $ with its derivatives \cite{plato2003concise,stoer2013,torres2020reduction}. Before continuing,  we need to consider the following multi-index notation. Let $\nset{N}_0$ be the set $\nset{N}\cup\set{0}$, if $\gamma \in \nset{N}_0^n$, then

\begin{eqnarray}
\left\{
\begin{array}{l}
 \gamma!:= \ds\prod_{k=1}^n [\gamma]_k ! \vspace{0.1cm}\\
 \abs{\gamma}:= \ds \sum_{k=1}^n [\gamma]_k\vspace{0.1cm}\\
 x^\gamma:= \ds \prod_{k=1}^n [x]_k^{[\gamma]_k}\vspace{0.1cm}\\
\der{\partial}{x}{\gamma}:= \dfrac{\partial^{\abs{\gamma}}}{\partial [x]_1^{[\gamma]_1}\partial [x]_2^{[\gamma]_2}\cdots \partial [x]_n^{[\gamma]_n} }
\end{array}\right. .
\end{eqnarray}

\begin{theorem}\label{teo:c2.01}
Let $ \Phi: \Omega \subset \nset{R}^n \to \nset {R}^n $ be an iteration function with a fixed point $ \xi \in \Omega $. Assuming that $\Phi $ is $ p$-times differentiable in $ \xi $ for some $ p \in \nset{N} $, and furthermore

\begin{eqnarray}\label{eq:c2.09}
\left\{
\begin{array}{cc}
\ds  \dfrac{\partial^\gamma [\Phi]_k(\xi) }{ \partial x^\gamma}=0, \ \forall k\geq 1 \mbox{ and } \forall  \abs{\gamma}<p, & \mbox{if }p\geq 2 \vspace{0.1cm}\\
\ds \norm{\Phi^{(1)}(\xi)}<1, & \mbox{if }p=1
\end{array}\right.,
\end{eqnarray}

where $\Phi^{(1)}$ denotes the \textbf{Jacobian matrix} of the function $\Phi$, then $ \Phi $ is (locally) convergent of (at least) order $ p $.

\begin{proof}

Let $\Phi:\nset{R}^n \to \nset{R}^n$ be an iteration function, and let $\set{\hat{e}_k}_{k=1}^n$ be the canonical basis of $\nset{R}^n$. Considering the following index notation (Einstein notation)

\begin{eqnarray*}
\Phi(x)=\sum_{k=1}^n [\Phi]_k(x)\hat{e}_k: = [\Phi]_k(x)\hat{e}_k=\hat{e}_k[\Phi]_k(x),
\end{eqnarray*}

and using the Taylor series expansion of a vector-valued function in multi-index notation, we obtain two cases:

\begin{itemize}

\item[i)]  Case $p\geq 2:$

\begin{align*}
\Phi(x_i)
=& \ds  \Phi(\xi)+  \sum_{\abs{\gamma} =1}^p \dfrac{1}{\gamma !}\hat{e}_k\dfrac{\partial^\gamma [\Phi]_k(\xi) }{ \partial x^\gamma}   (x_i-\xi)^\gamma   + \hat{e}_k[o]_k\left(\max_{\abs{\gamma}=p}  \set{ (x_i-\xi)^\gamma }\right) \\
=& \ds \Phi(\xi)+ \sum_{m =1}^p   \left( \sum_{\abs{\gamma} =m}\dfrac{1}{\gamma !} \hat{e}_k \dfrac{\partial^\gamma [\Phi]_k(\xi) }{ \partial x^\gamma}   (x_i-\xi)^\gamma \right)   + \hat{e}_k[o]_k\left(\max_{\abs{\gamma}=p}  \set{ (x_i-\xi)^\gamma }\right),
\end{align*}

then

\begin{align*}
\norm{\Phi(x_i)-\Phi(\xi)}&\leq  \ds \sum_{m =1}^p   \left( \sum_{\abs{\gamma} =m}\dfrac{1}{\gamma !}\norm{ \hat{e}_k \dfrac{\partial^\gamma [\Phi]_k(\xi) }{ \partial x^\gamma}   (x_i-\xi)^\gamma }\right)  +  \norm{\hat{e}_k  [o]_k\left(\max_{\abs{\gamma}=p}  \set{ (x_i-\xi)^\gamma }\right)} \\
&\leq   \ds \sum_{m =1}^p   \left( \sum_{\abs{\gamma} =m}\dfrac{1}{\gamma !}\norm{ \dfrac{\partial^\gamma [\Phi]_k(\xi) }{ \partial x^\gamma}  \hat{e}_k   }\right)\norm{x_i-\xi}^m+  o\left( \norm{x_i-\xi}^p \right), 
\end{align*}

assuming that $\xi$ is a fixed point of $\Phi$ and that $\dfrac{\partial^\gamma [\Phi]_k(\xi) }{ \partial x^\gamma}=0 \ \forall k\geq 1$ and $  \forall \abs{\gamma}<p$ is fulfilled, the previous expression implies that

\begin{eqnarray*}
\dfrac{\norm{\Phi(x_i)-\Phi(\xi)}}{\norm{x_i-\xi}^p}=\dfrac{\norm{x_{i+1}-\xi}}{\norm{x_i-\xi}^p}\leq\sum_{\abs{\gamma} =p}\dfrac{1}{\gamma !}\norm{ \dfrac{\partial^\gamma [\Phi]_k(\xi) }{ \partial x^\gamma}  \hat{e}_k   } +\dfrac{o\left(\norm{x_i-\xi}^p \right)}{\norm{x_i-\xi}^p},
\end{eqnarray*}

therefore

\begin{eqnarray*}
\lim_{i\to \infty} \dfrac{\norm{x_{i+1}-\xi}}{\norm{x_i-\xi}^p}\leq  \sum_{\abs{\gamma} =p}\dfrac{1}{\gamma !}\norm{ \dfrac{\partial^\gamma [\Phi]_k(\xi) }{ \partial x^\gamma}  \hat{e}_k   },
\end{eqnarray*}

as a consequence, if the sequence $\set{x_i}_{i=0}^\infty$ generated by \eqref{eq:2-001} converges to $\xi$, there exists a value $k>0$ such that

\begin{eqnarray*}
\norm{x_{i+1}-\xi}\leq \left( \sum_{\abs{\gamma} =p}\dfrac{1}{\gamma !}\norm{ \dfrac{\partial^\gamma [\Phi]_k(\xi) }{ \partial x^\gamma}  \hat{e}_k   }\right)\norm{x_i-\xi}^p, & \forall i\geq k,
\end{eqnarray*}

then $ \Phi $ is (locally) convergent of (at least) order $ p $.

 \item[ii)] Case $p=1:$

\begin{align*}
\Phi(x_i)=& \ds \Phi(\xi)+  \sum_{\abs{\gamma} =1}\dfrac{1}{\gamma !} \hat{e}_k \dfrac{\partial^\gamma [\Phi]_k(\xi) }{ \partial x^\gamma}   (x_i-\xi)^\gamma    + \hat{e}_k[o]_k\left(\max_{\abs{\gamma}=1}  \set{ (x_i-\xi)^\gamma }\right)\\
=&\Phi(\xi)+\Phi^{(1)}(x_i)(x_i-\xi)+\hat{e}_k[o]_k\left(\max_{\abs{\gamma}=1}  \set{ (x_i-\xi)^\gamma }\right),
\end{align*}

then

\begin{align*}
\norm{\Phi(x_i)-\Phi(\xi)}\leq  \norm{\Phi^{(1)}(\xi)} \norm{x_i-\xi }+   o\left( \norm{x_i-\xi} \right) ,
\end{align*}

assuming that $\xi$ is a fixed point of $\Phi$, the previous expression implies that

\begin{eqnarray*}
\dfrac{\norm{\Phi(x_i)-\Phi(\xi)}}{\norm{x_i-\xi}}=\dfrac{\norm{x_{i+1}-\xi}}{\norm{x_i-\xi}}\leq \norm{\Phi^{(1)}(\xi)} +\dfrac{o\left(\norm{x_i-\xi} \right)}{\norm{x_i-\xi}},
\end{eqnarray*}

therefore

\begin{eqnarray*}
\lim_{i\to \infty} \dfrac{\norm{x_{i+1}-\xi}}{\norm{x_i-\xi}}\leq\norm{\Phi^{(1)}(\xi)},
\end{eqnarray*}

as a consequence, if the sequence $\set{x_i}_{i=0}^\infty$ generated by \eqref{eq:2-001} converges to $\xi$, there exists a value $k>0$ such that

\begin{eqnarray*}
\norm{x_{i+1}-\xi}\leq \norm{\Phi^{(1)}(\xi)}\norm{x_i-\xi}, & \forall i\geq k,
\end{eqnarray*}

considering $m\geq 1$, from the previous inequality we obtain that

\begin{align*}
\norm{x_{i+m}-\xi}\leq &\norm{\Phi^{(1)}(\xi)}\norm{x_{i+m-1}-\xi}\leq \norm{\Phi^{(1)}(\xi)}^2\norm{x_{i+m-2}-\xi} \leq \cdots \leq\norm{\Phi^{(1)}(\xi)}^{m}\norm{x_{i}-\xi},
\end{align*}

and assuming that $\norm{\Phi^{(1)}(\xi)}<1$ is fulfilled

\begin{eqnarray*}
\lim_{m \to \infty}\norm{x_{i+m}-\xi}\leq \lim_{m \to \infty}\norm{\Phi^{(1)}(\xi)}^{m}\norm{x_{i}-\xi} \to 0,
\end{eqnarray*}

then $ \Phi $ is (locally) convergent of order (at least) linear.

\end{itemize}

\end{proof}

\end{theorem}

The next corollary follows from the previous theorem

\begin{corollary}\label{cor:2-001}
Let $\Phi:\nset{R}^n \to \nset{R}^n$ be an iteration function. If $\Phi$ defines a sequence $\set{x_i}_{i=0}^\infty$ such that $x_i\to \xi$, and if the following condition is true

\begin{eqnarray}\label{eq:c2.16}
\lim_{x\to \xi}\norm{\Phi^{(1)}(x)}\neq 0,
\end{eqnarray}

then $\Phi$ has an order of convergence (at least) linear  in $B(\xi;\delta)$.
\end{corollary}

\section{Riemann-Liouville Fractional Derivative}

One of the key pieces in the study of fractional calculus is the iterated integral, which is defined as follows \cite{hilfer00}

\begin{definition}
Let $ L_{loc} ^ 1 (a, b) $ be the space of locally integrable functions in the interval $ (a, b) $. If $ f $ is a function such that $ f \in L_ {loc} ^ 1 (a, \infty) $, then the $n$-th iterated integral of the function $ f $ is given by 

\begin{eqnarray}\label{eq:c1.16}
\begin{array}{c}
\ds \ifr{}{a}{I}{x}{n} f(x)=\ifr{}{a}{I}{x}{}\left(\ifr{}{a}{I}{x}{n-1} f(x)  \right)=\frac{1}{(n-1)!}\int_a^x(x-t)^{n-1}f(t)dt,
\end{array}
\end{eqnarray}

where

\begin{eqnarray*}
\ifr{}{a}{I}{x}{} f(x):=\int_a^x f(t)dt.
\end{eqnarray*}

\end{definition}

Considerate that $ (n-1)! = \gam{n} $
, a generalization of \eqref{eq:c1.16} may be obtained for an arbitrary order $ \alpha> 0 $

\begin{eqnarray}\label{eq:c1.17}
\ifr{}{a}{I}{x}{\alpha} f(x)=\dfrac{1}{\gam{\alpha}}\int_a^x(x-t)^{\alpha-1}f(t)dt,
\end{eqnarray}

the equation \eqref{eq:c1.17}  correspond to the definition of \textbf{Riemann-Liouville (right) fractional integral}.  Fractional integrals fulfill the \textbf{semigroup property}, which is given in the following proposition \cite{hilfer00}:

\begin{proposition}
Let $ f $ be a function. If $ f \in L_{loc} ^ 1 (a, \infty) $, then the fractional integrals of $ f $ fulfill that

\begin{eqnarray}\label{eq:c1.19}
\ifr{}{a}{I}{x}{\alpha} \ifr{}{a}{I}{x}{\beta}f(x) = \ifr{}{a}{I}{x}{\alpha + \beta}f(x),& \alpha,\beta>0.
\end{eqnarray}

\end{proposition}

From the previous proposition, and considering that the operator $ d / dx $  is the inverse operator to the left of the operator $ \ifr {}{a}{I}{x}{} $, any integral $ \alpha$-th of a function $ f \in L_{loc} ^ 1 (a, \infty) $ may be written as

\begin{eqnarray}\label{eq:c1.20}
\ifr{}{a}{I}{x}{\alpha}f(x)=\dfrac{d^n}{dx^n}\ifr{}{a}{I}{x}{n}\left( \ifr{}{a}{I}{x}{\alpha}f(x) \right)=\dfrac{d^n}{dx^n}\left( \ifr{}{a}{I}{x}{n+\alpha}f(x)\right).
\end{eqnarray}

With the previous results, we can build the operator  \textbf{Riemann-Liouville fractional derivative} as follows \cite{hilfer00,kilbas2006theory}

\begin{eqnarray}\label{eq:c1.23}
\normalsize
\begin{array}{c}
\ifr{}{a}{D}{x}{\alpha}f(x) := \left\{
\begin{array}{cc}
\ds \ifr{}{a}{I}{x}{-\alpha}f(x), &\mbox{if }\alpha<0 \vspace{0.1cm}\\  
\ds \dfrac{d^n}{dx^n}\left( \ifr{}{a}{I}{x}{n-\alpha}f(x)\right), & \mbox{if }\alpha\geq 0
\end{array}
\right.
\end{array}, 
\end{eqnarray}

where  $ n = \lceil \alpha \rceil $.

\subsection{Examples of the Riemann-Liouville Fractional Derivative}

Before continuing, it is necessary to define the Beta function and the incomplete Beta function \cite{arfken85}, which are defined as follows

\begin{eqnarray}
B(p,q):=\int_0^1 t^{p-1}(1-t)^{q-1}dt, & \ds B_r(p,q):=\int_0^r t^{p-1}(1-t)^{q-1}dt, 
\end{eqnarray}

where $p$ and $q$ are positive values. Considering the following proposition:

\begin{proposition}\label{prop:01}
Let $f$ be a function, with 

\begin{eqnarray*}
f(x)=(x-c)^\mu , & \mu>-1, & c\in \nset{R},
\end{eqnarray*}

then for all $\alpha \in \nset{R}\setminus \nset{Z}$, the Riemann-Liouville fractional derivative of the above function may be written as

\begin{eqnarray}
\normalsize
\begin{array}{c}
\ifr{}{a}{D}{x}{\alpha}f(x) = \left\{
\begin{array}{cc}
\ds \dfrac{\gam{\mu +1}}{\gam{\mu-\alpha+1}}(x-c)^{\mu-\alpha}G_{-\alpha}\left(\dfrac{a-c}{x-c},\mu+1\right), &\mbox{if }\alpha<0 \vspace{0.1cm}\\  
\ds \sum_{k=0}^n \binom{n}{k}  \dfrac{\gam{\mu +1}}{\gam{\mu +n-\alpha-k+1}}  (x-c)^{\mu +n-\alpha-k} G_{n-\alpha}^{(n-k)}\left(\dfrac{a-c}{x-c},\mu +1 \right), & \mbox{if }\alpha\geq 0
\end{array}
\right.
\end{array}, 
\end{eqnarray}

where

\begin{eqnarray}\label{eq:S2-1-002}
G_{\alpha}\left(\dfrac{a-c}{x-c},\mu+1 \right):= 1  - \dfrac{B_{\frac{a-c}{x-c}}(\mu+1,\alpha)}{B(\mu+1,\alpha)}.
\end{eqnarray}

\begin{proof}
The Riemann-Liouville fractional derivative of the function $f(x)$, through the equation \eqref{eq:c1.23}, presents two cases:

\begin{itemize}
\item[i)] If $\alpha <0$, then :

\begin{eqnarray*}
\ifr{}{a}{D}{x}{\alpha}f(x)=\dfrac{1}{\gam{-\alpha}}\int_a^x(x-t)^{-\alpha-1}(t-c)^\mu dt,
\end{eqnarray*}

taking the change of variable $t=c+(x-c)u$ in the previous expression

\begin{align*}
\ifr{}{a}{D}{x}{\alpha}f(x)=\dfrac{(x-c)^{\mu-\alpha}}{\gam{-\alpha}}\int_{\frac{a-c}{x-c}}^1(1-u)^{-\alpha-1}u^\mu du,
\end{align*}

the above result may be rewritten in terms of  the Beta function and the incomplete Beta function as follows

\begin{align*}
\ifr{}{a}{D}{x}{\alpha}f(x)
=&\dfrac{(x-c)^{\mu-\alpha}}{\gam{-\alpha}}\left( B(\mu+1,-\alpha)  -B_{\frac{a-c}{x-c}}(\mu+1,-\alpha) \right) \\
=& B(\mu +1,-\alpha)\dfrac{(x-c)^{\mu-\alpha}}{\gam{-\alpha}}\left( 1  - \dfrac{B_{\frac{a-c}{x-c}}(\mu+1,-\alpha)}{B(\mu+1,-\alpha)} \right),
\end{align*}

and considering \eqref{eq:S2-1-002}, we obtain that

\begin{eqnarray}\label{eq:S2-1-001}
\ifr{}{a}{D}{x}{\alpha}(x-c)^\mu= \dfrac{\gam{\mu +1}}{\gam{\mu-\alpha+1}}(x-c)^{\mu-\alpha}G_{-\alpha}\left(\dfrac{a-c}{x-c},\mu+1\right).
\end{eqnarray}

\item[ii)] If $\alpha \geq 0$, then:

\begin{eqnarray*}
\ifr{}{a}{D}{x}{\alpha}f(x)=\dfrac{1}{\gam{n-\alpha}}\dfrac{d^n}{dx^n}\int_a^x(x-t)^{n-\alpha-1}(t-c)^\mu dt,
\end{eqnarray*}

taking the change of variable $t=c+(x-c)u$ in the previous expression

\begin{align*}
\ifr{}{a}{D}{x}{\alpha}f(x)=\dfrac{1}{\gam{n-\alpha}}\dfrac{d^n}{dx^n} \left[ (x-c)^{\mu +n-\alpha} \int_{\frac{a-c}{x-c}}^1(1-u)^{n-\alpha-1}u^\mu du \right] ,
\end{align*}

the above result may be rewritten in terms of  the Beta function and the incomplete Beta function as follows

\begin{align*}
\ifr{}{a}{D}{x}{\alpha}f(x)
=&\dfrac{1}{\gam{n-\alpha}}\dfrac{d^n}{dx^n} \left[ (x-c)^{\mu+n-\alpha}\left( B(\mu +1,n-\alpha)  -B_{\frac{a-c}{x-c}}(\mu +1,n-\alpha)\right) \right] \\
=&\dfrac{B(\mu +1,n-\alpha)}{\gam{n-\alpha}}\dfrac{d^n}{dx^n} \left[ (x-c)^{\mu +n-\alpha}\left( 1  - \dfrac{B_{\frac{a-c}{x-c}}(\mu +1,n-\alpha)}{B(\mu +1,n-\alpha)} \right) \right] ,
\end{align*}

and considering \eqref{eq:S2-1-002}, we obtain that

\begin{align*}
\ifr{}{a}{D}{x}{\alpha}f(x)=&\dfrac{\gam{\mu +1}}{\gam{\mu +n-\alpha+1}}\dfrac{d^n}{dx^n} \left[ (x-c)^{\mu +n-\alpha}G_{n-\alpha}\left(\dfrac{a-c}{x-c},\mu +1 \right) \right] \\
=&\dfrac{\gam{\mu +1}}{\gam{\mu +n-\alpha+1}}  \sum_{k=0}^n \binom{n}{k}\left(\dfrac{d^k}{dx^k} (x-c)^{\mu +n-\alpha} \right) G_{n-\alpha}^{(n-k)}\left(\dfrac{a-c}{x-c},\mu +1 \right),
\end{align*}

taking into account that in the classical calculus

\begin{eqnarray*}
\dfrac{d^k}{dx^k}(x-c)^\mu =\dfrac{\mu !}{(\mu -k)!}(x-c)^{\mu-k}=\dfrac{\gam{\mu +1}}{\gam{\mu -k+1}}(x-c)^{\mu -k},
\end{eqnarray*}

therefore

\begin{align}\label{eq:S2-1-003}
\ifr{}{a}{D}{x}{\alpha}(x-c)^\mu
= \sum_{k=0}^n \binom{n}{k}  \dfrac{\gam{\mu +1}}{\gam{\mu +n-\alpha-k+1}}  (x-c)^{\mu +n-\alpha-k} G_{n-\alpha}^{(n-k)}\left(\dfrac{a-c}{x-c},\mu +1 \right).
\end{align}

\end{itemize}
\end{proof}

\end{proposition}

From the previous proposition, we can note that the Riemann-Liouville fractional derivative presents an explicit dependence of the value $n=\lceil \alpha \rceil$. However, there exists a particular case in which this dependence disappears, as shown in the following proposition:

\begin{proposition}\label{prop:02}
Let $f$ be a function, with

\begin{eqnarray*}
f(x)=(x-a)^\mu , & \mu>-1, & a\in \nset{R},
\end{eqnarray*}

then for all $\alpha\in \nset{R}\setminus \nset{Z}$, the Riemann-Liouville fractional derivative of the above function may be written in general form as

\begin{eqnarray}\label{eq:S2-1-004}
\ifr{}{a}{D}{x}{\alpha}(x-a)^\mu=\dfrac{\gam{\mu+1}}{\gam{\mu-\alpha+x}}(x-a)^{\mu-\alpha}.
\end{eqnarray}

\begin{proof}
To prove the validity of the previous equation for all $\alpha\in \nset{R}\setminus \nset{Z}$, it is necessary to note that from the \textbf{Proposition \ref{prop:01}} , the following limits may be obtained

\begin{eqnarray*}
\ifr{}{a}{D}{x}{\alpha}(x-a)^\mu=\lim_{c\to a}\ifr{}{a}{D}{x}{\alpha}(x-c)^\mu,
\end{eqnarray*}

\begin{eqnarray*}
\lim_{c \to a}G_{\alpha}\left(\dfrac{a-c}{x-c},m+1 \right)=G_\alpha(0,\mu+1)=1,
\end{eqnarray*}

then consider two cases:

\begin{itemize}
\item[i)] If  $\alpha < 0$, from the equation \eqref{eq:S2-1-001}, we obtain that

\begin{align*}
\ifr{}{a}{D}{x}{\alpha}(x-a)^\mu=&\dfrac{\gam{\mu +1}}{\gam{\mu-\alpha+1}} \lim_{c\to a} \left( (x-c)^{\mu-\alpha}G_{-\alpha}\left(\dfrac{a-c}{x-c},\mu+1\right) \right) \\
=&\dfrac{\gam{\mu +1}}{\gam{\mu-\alpha+1}} (x-a)^{\mu-\alpha}G_{-\alpha}\left(0,\mu+1\right) \\
=&\dfrac{\gam{\mu +1}}{\gam{\mu-\alpha+1}} (x-a)^{\mu-\alpha}.
\end{align*}

\item[i)] If  $\alpha \geq 0$, from the equation \eqref{eq:S2-1-003}, we obtain that

\begin{align*}
\ifr{}{a}{D}{x}{\alpha}(x-a)^\mu=&\sum_{k=0}^n \binom{n}{k}  \dfrac{\gam{\mu +1}}{\gam{\mu +n-\alpha-k+1}}\lim_{c\to a}\left(   (x-c)^{\mu +n-\alpha-k} G_{n-\alpha}^{(n-k)}\left(\dfrac{a-c}{x-c},\mu +1 \right) \right) \\
=&\sum_{k=0}^n \binom{n}{k}  \dfrac{\gam{\mu +1}}{\gam{\mu +n-\alpha-k+1}}  (x-a)^{\mu +n-\alpha-k} G_{n-\alpha}^{(n-k)}\left(0,\mu +1 \right)  \\
=& \binom{n}{n}  \dfrac{\gam{\mu +1}}{\gam{\mu +n-\alpha-n+1}}  (x-a)^{\mu +n-\alpha-n} G_{n-\alpha}^{(0)}\left(0,\mu +1 \right)  \\
=&\dfrac{\gam{\mu +1}}{\gam{\mu-\alpha+1}} (x-a)^{\mu-\alpha}.
\end{align*}

\end{itemize}

\end{proof}

\end{proposition}

From the previous proposition, the following corollary is obtained

\begin{corollary}
Let $f:\Omega \subset \nset{R} \to \nset{R}$ be a function, with $f \in L_{loc}^1(a,\infty)$. Assuming furthermore that $f\in C^\infty(a,\infty)$,   such that $f$ may be written in terms of its Taylor series around the point $x=a$, that is,

\begin{eqnarray*}
f(x)=\sum_{k=0}^\infty \dfrac{f^{(k)}(a)}{k!}(x-a)^k,
\end{eqnarray*}

then for all $\alpha \in \nset{R}\setminus \nset{Z}$, the Riemann-Liouville fractional derivative of the aforementioned function, may be written as follows

\begin{eqnarray}
\ifr{}{a}{D}{x}{\alpha}f(x)=\sum_{k=0}^\infty \dfrac{f^{(k)}(a)}{\gam{k-\alpha+1}} (x-a)^{k-\alpha}.
\end{eqnarray}

\end{corollary}

Finally, applying the  operator \eqref{eq:c1.23} with $a=0$ to the  function $ x^{\mu} $, with $\mu> -1$, from the \textbf{Proposition \ref{prop:02}} we obtain the following result

\begin{eqnarray}\label{eq:c1.13}
\ifr{}{0}{D}{x}{\alpha}x^\mu = 
 \dfrac{\gam{\mu+1}}{\gam{\mu-\alpha+1}}x^{\mu-\alpha}, & \alpha\in \nset{R}\setminus \nset{Z}.
\end{eqnarray}

\section{Fractional Pseudo-Newton Method}

Let $f:\Omega \subset \nset{R}^n\to \nset{R}^n$ be a function. We can consider the problem of finding a value $\xi\in \Omega$ such that $\norm{f(\xi)}$=0. A first approximation to value $\xi$ is by a linear approximation of the function $f$ in a valor $x_i\in \Omega$ with $\norm{x_i- \xi} < \epsilon$, that is,

\begin{eqnarray}\label{eq:c2.36}
f(x)\approx f(x_i)+f^{(1)}(x_i)(x-x_i),
\end{eqnarray}

then considering that $ \xi $ is a zero of $ f $, from the previous expression we obtain that

\begin{eqnarray*}
0\approx f(x_i)+f^{(1)}(x_i)(\xi-x_i) & \Rightarrow &  \xi \approx x_i- \left(f^{(1)}(x_i) \right)^{-1} f(x_i),
\end{eqnarray*}

consequently, may be generated a sequence $ \set{x_i}_{i = 0} ^ \infty $ that approximates the value $ \xi $  using the iterative method

\begin{eqnarray*}
x_{i+1}:=\Phi(x_i)=x_i- \left(f^{(1)}(x_i) \right)^{-1} f(x_i), & i=0,1,2,\cdots,
\end{eqnarray*}

which corresponds to well-known Newton's method. However, the equation \eqref{eq:c2.36} is not the only way to generate a linear approximation to the  function $ f $ in the point $ x_i $, another alternative is to use the next approximation

\begin{eqnarray}\label{eq:c2.37}
f(x)\approx f(x_i)+mI_n(x-x_i),
\end{eqnarray}

where $ I_n $ corresponds to the identity matrix of $n\times n$ and $ m $ is any constant value of a slope, that allows the approximation  \eqref{eq:c2.37} to the  function $ f $ to be valid. The previous equation allows to obtain the following iterative method

\begin{eqnarray}\label{eq:c2.39}
x_{i+1}:=\Phi(x_i)= x_i- \left( m^{-1}I_n \right) f(x_i), & i=0,1,2\cdots,
\end{eqnarray}

which corresponds to a particular case of the \textbf{parallel chord method} \cite{ortega1970iterative}. It is necessary to mention that for some definitions of fractional derivative, it is fulfilled that the derivative of the order $ \alpha $ of a constant is different from zero, that is,

\begin{eqnarray}\label{eq:c2.30}
\partial_k^\alpha c :=\der{\partial}{[x]_k}{\alpha}c \neq 0 , & c=constant,
\end{eqnarray}

where $ \partial_k ^ \alpha $ denotes any fractional derivative applied only in the component $ k $, that does not cancel the constants (for example: Riesz, Grünwald–Letnikov, Riemann-Liouville, etc.\cite{miller93,hilfer00,oldham74,kilbas2006theory,brambila2017fractal, martinez2017applications1,martinez2017applications2}), and that fulfills the following continuity relation with respect to the order $ \alpha $ of the derivative

\begin{eqnarray}\label{eq:c2.301}
\lim_{\alpha \to 1}\partial_k^\alpha c=\partial_kc.
\end{eqnarray}

Considering a function $\Phi:(\nset{R}\setminus \nset{Z})\times \nset{C}^n \to \nset{C}^n$. Then, using as a basis the idea of the method \eqref{eq:c2.39}, and considering any fractional derivative that fulfills the conditions \eqref{eq:c2.30} and \eqref{eq:c2.301}, we can define the \textbf{fractional pseudo-Newton method} \cite{torres2020fractional,torres2020nonlinear}  as follows

\begin{eqnarray}\label{eq:c2.401}
x_{i+1}:=\Phi(\alpha, x_i)= x_i- P_{\epsilon,\beta}(x_i) f(x_i), & i=0,1,2\cdots,
\end{eqnarray}

in particular is possible to take $\alpha\in[-2,2]\setminus\nset{Z}$ \cite{torreshern2020}, where $ P_{\epsilon, \beta} (x_i) $ is a matrix evaluated in the value $ x_i $, which is given by the following expression

\begin{eqnarray}\label{eq:c2.402}
P_{\epsilon,\beta}(x_i):=\left([P_{\epsilon,\beta}]_{jk}(x_i)\right)=\left( \partial_k^{\beta(\alpha,[x_i]_k)}\delta_{jk}+ \epsilon\delta_{jk}  \right)_{x_i},
\end{eqnarray}

where

\begin{eqnarray}\label{eq:c2.403}
\partial_k^{\beta(\alpha,[x_i]_k)}\delta_{jk}:= \der{\partial}{[x]_k}{\beta(\alpha,[x_i]_k)}\delta_{jk}, & 1\leq j,k\leq n,
\end{eqnarray}

with $ \delta_{jk} $ the Kronecker delta, $ \epsilon $ a positive constant $ \ll 1 $, and $ \beta (\alpha, [x_i]_k) $ a function defined as follows

\begin{eqnarray}\label{eq:c2.34}
\beta(\alpha,[x_i]_k):=\left\{
\begin{array}{cc}
\alpha, &\mbox{if \hspace{0.1cm} } \sqrt{  [x_i]_k [\overline{x_i}]_k   }> 0 \vspace{0.1cm}\\
1,& \mbox{if \hspace{0.1cm}  }  \sqrt{  [x_i]_k [\overline{x_i}]_k   }=0
\end{array}\right..
\end{eqnarray}

Due to the part of the integral operator that fractional derivatives usually have, we consider in the matrix \eqref{eq:c2.402} that each fractional derivative is obtained for a real variable $[x]_k$, and if the result allows it, this variable is subsequently made to tend to a complex variable $[x_i]_k$, that is,

\begin{eqnarray}
P_{\epsilon,\beta}(x_i):=P_{\epsilon,\beta}(x)\bigg{|}_{ x\longrightarrow x_i}  , & x\in \nset{R}^n, & x_i\in \nset{C}^n.
\end{eqnarray}

It should be mentioned that the value $ \alpha = 1 $ in \eqref{eq:c2.34}, is taken to avoid the discontinuity that is generated when using the fractional derivative of constants in the value $ x = 0 $. Furthermore, since in the previous method  $\norm{\Phi^{(1)}(\alpha,\xi)}\neq 0$ if $\norm{f(\xi)}=0$, from the \textbf{Corollary \ref{cor:2-001}}, any sequence $ \set{x_i} _ {i = 0} ^ \infty $ generated by the iterative method \eqref {eq:c2.401} has an order of convergence (at least) linear.

To finish this section, it is necessary to mention that although the interest in fractional calculus has mainly focused on the study and development of techniques to solve differential equation systems of order non-integer \cite{miller93,hilfer00,oldham74,kilbas2006theory,brambila2017fractal, martinez2017applications1,martinez2017applications2}. Over the years, iterative methods have also been developed that use the properties of fractional derivatives to solve algebraic equation systems \cite{gao2009local,brambila2018fractional,gdawiec2019visual,gdawiec2020newton,akgul2019fractional,cordero2019variant,torreshern2020,
torres2020fractional,torres2020nonlinear}. These methods may be called \textbf{fractional iterative methods}, which under certain conditions, may accelerate their speed of convergence with the implementation of the Aitken's method     \cite{stoer2013,brambila2018fractional}. It should be noted that depending on the definition of fractional derivative used, fractional iterative methods have the particularity that they may be used of local form \cite{gao2009local} or  of global form \cite{torreshern2020}.

\subsection{Some Examples}

Instructions for implementing the method \eqref{eq:c2.401} along with information to provide values $\alpha \in [-2,2]\setminus \nset{Z}$ are found in the reference \cite{torreshern2020}. For rounding reasons, only for the examples the following function is defined

\begin{eqnarray}\label{eq:4-001}
\rnd{[x]_k}{m}:=\left\{
\begin{array}{cc}
\re{[x]_k},& \mbox{ if \hspace{0.1cm}} \abs{\im{[x]_k}}\leq 10^{-m}\vspace{0.1cm}\\
\left[x\right]_k,& \mbox{ if \hspace{0.1cm}} \abs{\im{[x]_k}}> 10^{-m}\vspace{0.1cm}\\
\end{array}\right..
\end{eqnarray}

Combining the function \eqref{eq:4-001} with the method \eqref{eq:c2.401}, the following iterative method is defined

\begin{eqnarray}\label{eq:c2.40}
x_{i+1}:=\rnd{\Phi(\alpha, x_i)}{5}, & i=0,1,2\cdots.
\end{eqnarray}

\begin{example}

Let $\set{f_k}_{k=0}^\infty$ be a sequence of functions, with

\begin{eqnarray*}
f_k(x)=-\gamma-\log\left(x\right)-\sum_{m=1}^k \dfrac{(-1)^m x^{2 m}}{2 m \Gamma(2 m+1)}      & \underset{k \to \infty}{\longrightarrow} & \int_x^\infty \dfrac{\cos(t)}{t}dt,
\end{eqnarray*}

where $\gamma$ is the Euler–Mascheroni constant \cite{arfken85}.  Then considering the value $k=50$, the initial condition $x_0=0.018$ is chosen to use the iterative method given by \eqref{eq:c2.40} along with fractional derivative given by \eqref{eq:c1.13}. Consequently, we obtain the results of the Table \ref{tab:04}

\begin{table}[!ht]
\centering
\footnotesize
$
\begin{array}{c|ccccc}
\toprule
&\alpha& x_n&\norm{x_n - x_{n-1} }_2  &\norm{f_{50}\left(x_n \right)}_2& n \\ 
\midrule
1	&	-1.11419	&	15.7703495        	&	6.60000e-7	&	9.53989e-9	&	33	\\
 2	&	 -1.10571	&	  22.03613996        	&	3.20000e-7	&	3.67622e-9	&	  29	\\
3 	&	 -1.09148	&	   3.38418043        	&	3.00000e-7	&	2.13662e-9	&	  23	\\
 4	&	 -1.09147	&	   9.52557552        	&	4.40000e-7	&	6.51952e-9	&	  23	\\
 5	&	 -1.00855	&	  28.30949579        	&	1.10000e-7	&	8.74154e-7	&	 110	\\
 6	&	  0.24869	&	 -18.65007933 + 4.7866359i 	&	5.12250e-7	&	9.98818e-7	&	  61	\\
 7	&	  0.25206	&	  -24.9715585 + 5.07028087i	&	4.21900e-7	&	8.51586e-7	&	  61	\\
 8	&	  0.25269	&	 -12.29980885 + 4.3899063i 	&	5.19230e-7	&	9.35313e-7	&	  61	\\
 9	&	  0.25298	&	 -18.65007948 - 4.78663644i	&	4.41022e-7	&	8.66794e-7	&	  64	\\
 10 	&	  0.28661	&	  -5.86092713 + 3.72436638i	&	5.74282e-7	&	9.63373e-7	&	  57	\\
11	&	  0.29031	&	 -24.97155803 - 5.07028097i	&	3.26497e-7	&	7.97508e-7	&	  64	\\
12	&	  0.30109	&	 -12.29980928 - 4.38990586i	&	4.44072e-7	&	9.93021e-7	&	  66	\\
13	&	  0.31605	&	  -5.86092745 - 3.72436648i	&	4.79270e-7	&	8.82199e-7	&	  62	\\
14	&	  1.11281	&	  31.44800139        	&	2.00000e-8	&	1.56267e-7	&	   3	\\
15	&	  1.14602	&	   0.61650487   	&	2.60000e-7	&	8.14734e-7	&	  33	\\
\bottomrule
\end{array}
$
\caption{Results obtained using the iterative method \eqref{eq:c2.40} with $\epsilon=e-3$.}\label{tab:04}
\end{table}

\end{example}

\begin{example}

Let $\set{f_k}_{k=0}^\infty$ be a sequence of functions, with

\begin{eqnarray*}
f_k(x)=\dfrac{\pi}{2}-\sum_{m=0}^k \dfrac{(-1)^m x^{2 m+1}}{(2 m+1) \Gamma(2 m+2)} & \underset{k \to \infty}{\longrightarrow} & \int_x^\infty \dfrac{\sin(t)}{t}dt.
\end{eqnarray*}

Then considering the value $k=50$, the initial condition $x_0=1.85$ is chosen to use the iterative method  given by \eqref{eq:c2.40} along with fractional derivative given by \eqref{eq:c1.13}. Consequently, we obtain the results of the Table \ref{tab:03}

\begin{table}[!ht]
\centering
\footnotesize
$
\begin{array}{c|ccccc}
\toprule
&\alpha& x_n&\norm{x_n - x_{n-1} }_2  &\norm{f_{50}\left(x_n \right)}_2& n \\ 
\midrule
1	&	-0.83718	&	23.60399266        	&	4.10000e-7	&	9.80551e-9	&	30	\\
 2	&	 -0.81526	&	  29.87824476        	&	7.30000e-7	&	8.19133e-7	&	 273	\\
3 	&	 -0.71339	&	  17.33566366        	&	4.10000e-7	&	2.48591e-8	&	  34	\\
 4	&	 -0.71324	&	  11.08303768        	&	3.70000e-7	&	2.67499e-8	&	  29	\\
 5	&	 -0.71174	&	   4.89383571        	&	5.10000e-7	&	4.87621e-8	&	  24	\\
 6	&	  0.36251	&	 -12.29964074 - 4.38965942i	&	4.41814e-7	&	9.42697e-7	&	  42	\\
 7	&	  0.36333	&	 -31.27978791 - 5.29112884i	&	4.31045e-7	&	5.52647e-7	&	 391	\\
 8	&	  0.36684	&	 -24.97153098 - 5.07020771i	&	3.56931e-7	&	7.71425e-7	&	  42	\\
 9	&	  0.36976	&	 -18.65001989 - 4.78651292i	&	3.80132e-7	&	8.88913e-7	&	  40	\\
 10	&	  0.38451	&	 -24.97153097 + 5.07020788i	&	2.05913e-7	&	5.36982e-7	&	  44	\\
 11 	&	  0.38646	&	 -18.65002028 + 4.78651268i	&	2.30000e-7	&	5.54454e-7	&	  42	\\
12	&	  0.40682	&	 -12.29964106 + 4.38965969i	&	4.49110e-7	&	9.97466e-7	&	  37	\\
13	&	  0.44711	&	  -5.86005858 - 3.72373544i	&	4.72017e-7	&	9.90371e-7	&	  41	\\
14	&	  0.48437	&	  -5.86005854 + 3.72373556i	&	4.52217e-7	&	9.82000e-7	&	  38	\\
15	&	  0.55885	&	 -31.27978639 + 5.2911368i 	&	5.50000e-7	&	3.69789e-7	&	 183	\\
16	&	  1.41172	&	   1.92644561  	&	1.10000e-7	&	9.97696e-7	&	 196	\\
\bottomrule
\end{array}
$
\caption{Results obtained using the iterative method \eqref{eq:c2.40} with $\epsilon=e-3$.}\label{tab:03}
\end{table}

\end{example}

\newpage

\begin{example}

Let $f$ be a function, with

\begin{eqnarray*}
f(x)=
\begin{pmatrix}
\dfrac{1}{2} [x]_1\big{(} \sin\left( [x]_1 [x]_2 \big{)} -1\right)-\dfrac{1}{4\pi}[x]_2 \vspace{0.1cm}\\
\left( 1- \dfrac{1}{4\pi} \right)\left(e^{2[x]_1}-e\right) +e\left(\dfrac{1}{\pi}[x]_2-2[x]_1\right)
\end{pmatrix}.
\end{eqnarray*}

Then the initial condition $x_0=(0.86,0.86)^T$ is chosen to use the iterative method  given by \eqref{eq:c2.40} along with fractional derivative given by \eqref{eq:c1.13}. Consequently, we obtain the results of the Table \ref{tab:06}

\begin{table}[!ht]
\centering
\footnotesize
$
\begin{array}{c|cccccc}
\toprule
&\alpha& [x_n]_1& [x_n]_2 &\norm{x_n - x_{n-1} }_2  &\norm{f\left(x_n \right)}_2& n \\ 
\midrule
1	&	0.69508	&	1.01828092 + 0.52158397i	&	5.18478971 - 3.76689418i	&	1.15758e-7	&	7.24108e-7	&	48	\\
 2	&	 0.69632	&	 -0.13780201 + 0.87180277i	&	   2.16460973 + 4.68221216i	&	1.31909e-7	&	9.55148e-7	&	  80	\\
3 	&	 0.7283 	&	 -0.13780202 - 0.87180273i	&	   2.16460988 - 4.68221226i	&	9.11043e-8	&	8.81449e-7	&	 100	\\
 4	&	 0.72889	&	 -0.15442216        	&	   1.14021866        	&	6.22977e-7	&	8.30511e-7	&	  60	\\
 5	&	 0.75757	&	  1.01828092 - 0.5215839i 	&	   5.18479004 + 3.76689413i	&	9.21954e-8	&	9.83408e-7	&	  80	\\
 6	&	 0.78188	&	 -0.20477864 - 1.30850366i	&	   2.21623485 - 7.86783099i	&	5.56776e-8	&	9.92736e-7	&	 246	\\
 7	&	 0.82863	&	  1.14584377 - 0.68994257i	&	   8.09450013 + 5.9960712i 	&	3.31662e-8	&	9.77549e-7	&	 193	\\
 8	&	 0.86097	&	  1.14584377 + 0.68994256i	&	   8.09450017 - 5.99607116i	&	2.64575e-8	&	9.42041e-7	&	 249	\\
 9	&	 1.11159	&	  1.70987637        	&	 -18.87534307        	&	1.41421e-8	&	9.92487e-7	&	 447	\\
 10	&	 1.14766	&	  1.48216448        	&	  -8.41311536        	&	1.41421e-8	&	8.86632e-7	&	 233	\\
 11 	&	 1.17262	&	 -1.36674692 + 0.07786741i	&	     -5.76423 + 0.47853094i	&	2.00000e-8	&	9.92337e-7	&	 394	\\
12	&	 1.18538	&	 -1.36674698 - 0.07786726i	&	  -5.76422966 - 0.4785315i 	&	2.23607e-8	&	9.88600e-7	&	 387	\\
13	&	 1.19954	&	  1.57643706        	&	   -12.098725        	&	1.41421e-8	&	7.09538e-7	&	 386	\\
14	&	 1.20058	&	  1.64946521        	&	 -15.55495398        	&	1.41421e-8	&	9.10544e-7	&	 465	\\
15	&	 1.2852 	&	 -0.76073057 + 0.14192444i	&	  -2.11123992 + 0.82667655i	&	1.02470e-7	&	8.39720e-7	&	  97	\\
16	&	 1.28668	&	 -0.76073047 - 0.14192446i	&	  -2.11123884 - 0.8266763i 	&	2.74044e-7	&	7.74152e-7	&	  81	\\
17	&	 1.29642	&	  1.34362303    	&	  -4.29400761      	&	7.61577e-8	&	4.60872e-7	&	  92	\\
\bottomrule
\end{array}
$
\caption{Results obtained using the iterative method \eqref{eq:c2.40} with $\epsilon=e-3$.}\label{tab:06}
\end{table}

\end{example}

\section{Approximation to zeros of the Riemann zeta function}

A systematic study of the Riemann zeta function is beyond the intent of this document. However, the basic information necessary to approximate their zeros  using  the iterative method \eqref{eq:c2.40} will be presented.    Let $\zeta:\Omega \subset \nset{C}\setminus \set{1}\to \nset{C} $ be the Riemann zeta function with $\Omega=\set{x\in \nset{C} \ : \ \re{x}>1}$. The function $\zeta$ is defined as follows

\begin{eqnarray}
\zeta(x)=\sum_{k=1}^\infty \dfrac{1}{k^x},
\end{eqnarray}

the previous expression may be extended for all $x\in \nset{C}\setminus \set{1}$  by \textbf{analytic continuation}. With which it is possible to obtain the following functional equation 

\begin{eqnarray}
\zeta(x)=2^x \pi^{x-1}\sin \left(\dfrac{\pi x}{2} \right)\Gamma(1-x)\zeta(1-x).
\end{eqnarray}

On the other hand, there exists a series version of the Riemann zeta function, which has the characteristic of being globally convergent for all $x\in \nset{C}\setminus \set{1}$. This version of the function $\zeta$ was conjectured by Konrad Knopp and proved by Helmut Hasse in 1930 \cite{joyner2011selected},  and it is given by the following expression

\begin{eqnarray}\label{eq:5-001}
\zeta (x)=\dfrac {1}{1-2^{1-x}}\sum _{m=0}^{\infty }{\frac {1}{2^{m+1}}}\sum _{p=0}^{m}(-1)^{p}{m \choose p}\left( p+1\right) ^{-x}.
\end{eqnarray}

It is necessary to mention that the expression \eqref{eq:5-001} is very useful, because it allows us to make numerical approximations of the nontrivial zeros of the Riemann zeta function.

\begin{example}

Let $\set{f_k}_{k=0}^\infty$ be a sequence of functions, with

\begin{eqnarray*}
f_k(x)=\dfrac {1}{1-2^{1-x}}\sum _{m=0}^{k }{\frac {1}{2^{m+1}}}\sum _{p=0}^{m}(-1)^{p}{m \choose p}\left(p+1\right)^{-x}    & \underset{k \to \infty}{\longrightarrow} & \zeta(x).
\end{eqnarray*}

Then considering the value $k=50$, the initial condition $x_0=0.5 + 31.51i$ is chosen to use the iterative method  given by \eqref{eq:c2.40} along with fractional derivative given by \eqref{eq:c1.13}. Consequently, we obtain the results of the Table \ref{tab:05}

\begin{table}[!ht]
\centering
\footnotesize
$
\begin{array}{c|ccccc}
\toprule
&\alpha& x_n&\norm{x_n - x_{n-1} }_2  &\norm{f_{50}\left(x_n \right)}_2& n \\ 
\midrule
    1     & -1.13348 & -5.99997526 & 9.90000e-7 & 1.46057e-7 & 266 \\
    2     & -0.94547 & -2.00001275 & 8.00000e-7 & 3.88212e-7 & 146 \\
    3     & -0.88567 & -9.9999519 & 7.90000e-7 & 5.72079e-7 & 54 \\
    4     & -0.15505 & 0.49999963 + 14.13472531i & 9.89798e-7 & 3.22385e-7 & 22 \\
    5     & -0.06835 & 0.54743578 + 947.15979693i & 8.59593e-7 & 3.79999e-7 & 349 \\
    6     & -0.06585 & 0.48286572 - 305.3825058i & 8.95377e-7 & 4.74515e-7 & 251 \\
    7     & -0.06135 & 0.47107264 - 430.95892398i & 9.80051e-7 & 4.94661e-7 & 108 \\
    8     & -0.02125 & 0.5460791 + 321.6875341i & 7.92023e-7 & 5.38613e-7 & 61 \\
    9     & -0.01965 & 0.48286548 + 305.38250571i & 6.32456e-7 & 1.92115e-7 & 49 \\
    10    & -0.01685 & 0.61664333 + 276.09845269i & 5.95483e-7 & 2.99694e-7 & 99 \\
    11    & -0.01275 & 0.58077376 + 224.41293459i & 9.49210e-7 & 8.04311e-7 & 267 \\
    12    & -0.01015 & 0.56816592 + 192.65292746i & 9.39628e-7 & 7.02867e-7 & 192 \\
    13    & -0.00755 & 0.46159425 + 158.41373653i & 8.31505e-7 & 4.66903e-7 & 175 \\
    14    & -0.00545 & 0.56284646 + 133.97913455i & 9.39840e-7 & 7.30271e-7 & 213 \\
    15    & -0.00465 & 0.48007233 + 123.75232554i & 9.48472e-7 & 8.25385e-7 & 278 \\
    16    & -0.00195 & 0.4913738 + 95.32674434i & 5.73847e-7 & 1.78959e-7 & 172 \\
    17    & 0.01795 & 0.50002627 - 37.5861995i & 9.40213e-7 & 9.17598e-7 & 255 \\
    18    & 0.02145 & 0.49137435 - 95.32674404i & 9.88130e-7 & 3.70864e-7 & 125 \\
    19    & 0.02165 & 0.50000116 - 32.93506385i & 8.62670e-7 & 6.14360e-7 & 120 \\
    20    & 0.03195 & 0.49999981 - 25.01085752i & 7.35391e-7 & 3.11793e-7 & 53 \\
    21    & 0.03295 & 0.46159416 - 158.41373657i & 8.10555e-7 & 3.57361e-7 & 139 \\
    22    & 0.03305 & 0.50000005 - 30.42487677i & 9.08020e-7 & 6.57518e-7 & 73 \\
    23    & 0.03425 & 0.50007248 - 40.91868687i & 6.35059e-7 & 3.03925e-7 & 87 \\
    24    & 0.03785 & 0.50000027 - 21.02203963i & 9.83107e-7 & 3.08350e-7 & 31 \\
    25    & 0.04495 & 0.50000036 - 14.13472527i & 8.24197e-7 & 3.03119e-7 & 21 \\
    26    & 0.06185 & 0.50000011 + 21.02203962i & 3.21403e-7 & 1.28471e-7 & 40 \\
    27    & 0.06735 & 0.40419552 + 1185.59789987i & 5.52268e-7 & 4.10670e-7 & 67 \\
    28    & 0.14095 & 0.50689898 - 2878.59861379i & 3.55106e-7 & 9.12771e-7 & 330 \\
    29    & 0.24136 & 0.44802022 - 495.84458602i & 2.32594e-7 & 8.59529e-7 & 60 \\
    30    & 0.24176 & 0.54285436 - 466.79078671i & 2.08806e-7 & 9.24091e-7 & 218 \\
    31    & 0.24256 & 0.58122569 - 404.58829656i & 2.62488e-7 & 9.89996e-7 & 354 \\
    32    & 0.24306 & 0.59895653 - 373.61214166i & 2.26716e-7 & 9.90256e-7 & 464 \\
    33    & 0.24346 & 0.46960335 - 346.68847181i & 2.10950e-7 & 8.74078e-7 & 496 \\
    34    & 0.24766 & 0.46263812 - 146.98942449i & 3.08869e-7 & 8.33117e-7 & 496 \\
    35    & 0.24806 & 0.56284605 - 133.97913398i & 3.61248e-7 & 9.02931e-7 & 478 \\
    36    & 0.24836 & 0.48007222 - 123.75232428i & 2.72029e-7 & 7.75292e-7 & 484 \\
    37    & 0.25106 & 0.50076291 - 56.43890631i & 3.64005e-7 & 8.55540e-7 & 354 \\
    38    & 0.25156 & 0.50021423 - 43.32696237i & 4.86004e-7 & 7.09618e-7 & 312 \\
    39    & 0.26886 & 0.50000042 + 25.01085761i & 3.53553e-7 & 5.54795e-7 & 94 \\
    40    & 0.28486 & 0.50007238 + 40.91868722i & 3.13847e-7 & 8.22684e-7 & 51 \\
    41    & 0.29196 & 0.50002634 + 37.5861994i & 3.74433e-7 & 9.61617e-7 & 46 \\
    42    & 0.29836 & 0.50021436 + 43.32696223i & 2.56125e-7 & 5.52253e-7 & 42 \\
    43    & 0.30096 & 0.50125709 + 49.77289402i & 2.41661e-7 & 6.80998e-7 & 63 \\
    44    & 0.30126 & 0.49834763 + 52.96633019i & 2.45153e-7 & 4.93711e-7 & 43 \\
    45    & 0.30136 & 0.5007636 + 56.43890622i & 2.72029e-7 & 9.06750e-7 & 69 \\
    46    & 0.30196 & 0.49914035 + 60.81054836i & 2.87924e-7 & 9.28287e-7 & 66 \\
    47    & 0.30686 & 0.49914001 - 60.81054769i & 2.86531e-7 & 9.54106e-7 & 59 \\
    48    & 0.30696 & 0.47366839 - 67.04250707i & 2.62488e-7 & 8.53872e-7 & 94 \\
    49    & 0.30776 & 0.49834757 - 52.96632988i & 2.64764e-7 & 5.86654e-7 & 288 \\
    50    & 0.30786 & 0.50125707 - 49.77289365i & 2.34094e-7 & 6.76399e-7 & 325 \\
    51    & 0.31156 & 0.47366808 + 67.04250713i & 2.78927e-7 & 9.54910e-7 & 44 \\
    52    & 0.32336 & 0.50000099 + 30.42487725i & 2.59422e-7 & 9.16200e-7 & 120 \\
    53    & 0.32386 & 0.50000112 + 32.93506408i & 2.73130e-7 & 7.59689e-7 & 40 \\
\bottomrule
\end{array}
$
\caption{Results obtained using the iterative method \eqref{eq:c2.40} with $\epsilon=e-3$.}\label{tab:05}
\end{table}

\end{example}

The comments below are made under the assumption that the Riemann hypothesis given in equation \eqref{eq:1-004} is true. Table \ref{tab:05} shows certain $x_n$ values with $\abs{\re{x_n}- 0.5}\geq 0.1$, this is partly a consequence of the approximation of the function $\zeta$ by means of the function $f_k$, because in general the following condition is true

\begin{eqnarray*}
\kr{f_k} \neq \kr{\zeta} & \Leftrightarrow & k<\infty.
\end{eqnarray*}

In addition, it is necessary to mention that the use of iterative methods does not guarantee that we can get as near as we want to the value $ \xi $, normally it is only possible to determine a value $ x_n $ near to the value $ \xi $ given by the following expression

\begin{eqnarray}\label{eq:6-001}
x_n=\xi+ \delta_\xi, & \delta_\xi=\delta_\xi(n), &\norm{\delta_\xi}<1,
\end{eqnarray}

considering \eqref{eq:6-001}, it is necessary to give a definition that allows us to characterize the behavior of a function with respect to the values $ x_n $ in $\overline{B}( \xi;\delta_\xi)$.

\begin{definition}
Let $f:\Omega \subset \nset{R}^n \to \nset{R}^n$ be a function with a value $ \xi \in \Omega$ such that $\norm{f(\xi)}=0$. If when doing $ \xi \to \xi +  \delta_\xi  $, with  $ \norm{\delta_\xi} <1 $, it holds that

\begin{eqnarray}\label{eq:6-002}
\norm{f(\xi+ \delta_\xi)}= \norm{\delta_f}<1,
\end{eqnarray}

then the function $f$ is \textbf{(locally) stable} with respect to the value $ \xi $ in $\overline{B}(\xi,\delta_\xi)$.

\end{definition}

The condition \eqref{eq:6-002}, implies that for a function $ f $ to be  (locally) stable, it is necessary that a slight perturbation $ \delta_ \xi $ in its zeros does not generate a great perturbation $ \delta_f $ in its images. To try to observe the stability of function $\zeta$, we may consider $\delta=10^{-12}$ and the following values

\begin{eqnarray*}
\begin{array}{c|c}
\begin{array}{ccc}
x_{n_1}=-40 -1\delta& \Rightarrow & \abs{\zeta(x_{n_1})}\approx 4.8\times 10^{3}\vspace{0.1cm} \\
x_{n_2}=-40 +0 \delta& \Rightarrow & \abs{\zeta(x_{n_2})}= 0\vspace{0.1cm} \\
x_{n_3}=-40+1\delta & \Rightarrow & \abs{\zeta(x_{n_3})}\approx 4.8\times 10^{3}
\end{array}&
\begin{array}{ccc}
x_{n_1}=-60 -1\delta& \Rightarrow & \abs{\zeta(x_{n_1})}\approx 5.3\times 10^{21}\vspace{0.1cm} \\
x_{n_2}=-60 +0 \delta& \Rightarrow & \abs{\zeta(x_{n_2})}= 0 \vspace{0.1cm} \\
x_{n_3}=-60+1\delta & \Rightarrow & \abs{\zeta(x_{n_3})}\approx 5.3\times 10^{21}
\end{array}
\end{array}.
\end{eqnarray*}

By performing multiple examples, it is possible to show for the trivial zeros of function $\zeta$ the validity of the following affirmation:

\begin{eqnarray*}
\mbox{ If }\xi=-2x \mbox{ with  }x\in \nset{N} & \Rightarrow & \zeta \mbox{ is (locally) unstable in }\xi \mbox{ if } x\to \infty.
\end{eqnarray*}

The previous affirmation implies that it is complicated to make numerical approximations to trivial zeros of $\zeta$ by iterative methods using an initial condition $x_0 \neq -2x$ with $x\in \nset{N}$ when $x\to \infty$. Furthermore, since the behavior of the nontrivial zeros of function $\zeta$ is not completely determined, it is possible to generate the following question:

\begin{eqnarray*}
\mbox{ If }\xi=0.5+xi \mbox{ with }x\in \nset{R} & \Rightarrow & \mbox{Is }\zeta \mbox{ (locally) unstable in }\xi \mbox{ if } \abs{x}\to \infty \mbox{?}
\end{eqnarray*}

\section{Conclusions}

Due that the derivatives of constant functions are identically zero functions in conventional calculus, it is difficult to imagine them as tools that may be used to approximate the zeros of more complicated functions. However, in the fractional calculus, the derivatives of constant functions may be used to approximate the zeros of more complicated functions and even non-differentiables. Furthermore, it is necessary to mention that constant functions are the easiest functions to work with in fractional calculus, since the fractional derivatives of more complicated functions are expressed in many occasions in terms of infinite series, given by Mittag-Leffler functions or hypergeometric functions. In this document, as far as the authors know, an approximation to the zeros of the Riemann zeta function has been obtained for the first time using only derivatives of constant functions, which was possible only because a fractional iterative method was used, specifically the fractional pseudo-Newton method.

The fractional iterative methods are efficient at finding multiple zeros of a function using a single initial condition, in addition to having the particularity of finding complex zeros of polynomials using real initial conditions. For this reason, they are suitable iterative methods for functions that have a large number of zeros, as is the case with the Riemann zeta function. It is necessary to mention that the fractional pseudo-Newton method does not explicitly depend on the fractional derivative of the function for which zeros are sought, then it is an ideal iterative method for working with functions that can be expressed in terms of a series or for solving nonlinear systems in several variables.

\bibliography{Biblio}
\bibliographystyle{unsrt}
\nocite{yang2018non}
\nocite{yang2019nontrivial}

\end{document}